\documentclass
[11pt]{article}
\usepackage{amssymb}
\usepackage{mathrsfs}
\usepackage{amsfonts}
\usepackage{amsmath}

\usepackage[top=1in, bottom=1in, left=1.25in, right=1.25in]{geometry}
\usepackage{graphicx}
\begin{document}

\title{\textbf{Opinion Dynamics and Influencing on Random Geometric Graphs}}
\author{Weituo Zhang$^{1,\ast}$, Chjan Lim$^{1}$, Gyorgy Korniss$^{2}$ and Boleslaw Szymanski$^{3}$}

\maketitle

\begin{flushleft}
$^{\bf{1}}$ Department of Mathematical Sciences, Rensselaer Polytechnic Institute,
110 8$^{th}$ Street, Troy, NY, 12180-3590 USA \\
$^{\bf{2}}$ Department of Physics, Applied Physics and Astronomy, Rensselaer Polytechnic Institute,
110 8$^{th}$ Street, Troy, NY, 12180-3590 USA \\
$^{\bf{3}}$ Department of Computer Science,
Rensselaer Polytechnic Institute, 110 8$^{th}$ Street, Troy, NY, 12180-3590 USA \\
$^{\ast}$ E-mail: zhangw14@rpi.edu
\end{flushleft}

\begin{abstract}
We investigate the two-word Naming Game on two-dimensional random
geometric graphs. Studying this model advances our understanding of the
spatial distribution and propagation of opinions in social dynamics.
A main feature of this model is the spontaneous emergence of spatial
structures called opinion domains which are geographic regions  with clear boundaries
within which all individuals share the same opinion. We provide the
mean-field equation for the underlying dynamics and discuss several properties
of the equation such as the stationary solutions and two-time-scale
separation. For the evolution of the opinion domains we find that
the opinion domain boundary propagates at a speed proportional to
its curvature. Finally we investigate the impact of
committed agents on opinion domains and find the scaling of
consensus time.
\end{abstract}


\noindent
Relevant features of social and opinion dynamics
\cite{Galam_IJMPC2008,Durlauf_PNAS1999,Castellano_RMP2009} can be
investigated by prototypical agent-based models such as the voter
model \cite{Liggett,Castellano}, the Naming Game
\cite{Baronchelli1,Dall'Asta2}, or the majority model
\cite{Galam_PhysA1999,Krapivsky_PRL2003}. These models typically
include a large number of individuals, each of which is assigned
a state defined by the social opinions that it accepts and updates its state by interacting
with its neighbors. Opinion dynamics driven by local communication
on geographically embedded networks is of great interest to
understanding the spatial distribution and propagation of opinions. In
this paper we investigate the Naming Game (NG) on random geometric
graphs as a minimum model of this type. We focus on the NG but will
also compare it with other models of opinion dynamics.

The NG \cite{Baronchelli1,Dall'Asta2} was originally introduced in the
context of linguistics and spontaneous emergence of shared
vocabulary among artificial agents \cite{Steels, Kirby} to
demonstrate how autonomous agents can achieve global agreement
through pair-wise communications without a central coordinator.
Here, we employ a special version of the NG, called the two-word
\cite{Castello_EPJB2009,Xie,Xie_PLOS2012,Zhang1} Listener-Only
Naming Game (LO-NG) \cite{Zhang1,Baronchelli_PRE2011}. In
this version of the NG, each agent can either adopt one of the two
different opinions A, B, or take the neutral stand represented by their union, AB.
In each communication, a pair of neighboring agents are randomly chosen, the first one as the speaker and the second one as the listener.
The speaker holding A or B opinion will transmit its own opinion and the neutral speaker will transmit either A or B opinion with
equal probability. The listener holding A or B opinion will become neutral when it hears an opinion different from its own and the neutral listener will adopt whatever it hears. Detailed instances are shown in Supplementary Table.

Consensus formation
in the original NG (and its variations) on various regular and
complex networks have been studied
\cite{Baronchelli1,Dall'Asta2,Xie,Baronchelli2,Dall'Asta1,Lu1,Lu2,Lu3,Zhang2}.
In particular, the spatial and temporal scaling properties have been
analyzed by direct simulations and scaling arguments in
spatially-embedded regular and random (RGG) graphs
\cite{Baronchelli2,Lu2,Lu3}. These results indicated
\cite{Baronchelli2,Lu2,Lu3} that the consensus formation in these
systems is analogous to coarsening \cite{Bray}. In this paper, we
further elucidate on the emerging coarsening dynamics in the
two-word LO-NG on RGG by developing mean-field (or coarse-grained)
equations for the evolution of opinions. Our method of relating
microscopic dynamics to macroscopic behavior shares similar features
with those leading to effective Fokker-Planck and Langevin equations
in a large class of opinion dynamic models (including generalized
voter models with intermediate states)
\cite{Vazquez,DallAsta_JPA2008}.

A random geometric graph (RGG), also referred to as a spatial Poisson or
Boolean graph, models spatial effects explicitly
and therefore is of both technological and intellectual importance
\cite{Penrose, Collier}. In this model, each node is randomly
assigned geographic coordinates and two nodes are connected
if the distance between them is within the interaction radius $r$.
Another type of network with geographic information is the regular
lattice. Fundamental models for opinion dynamics on regular lattice
has been intensively studied \cite{Castellano_RMP2009,Liggett,Cox}.
In many aspects, opinion
dynamics behaves similarly on RGGs and regular lattices with the
same dimensionality, but in our study, we also observe several
differences. For example, the length scale of spatial coarsening
for large $t$ is $l(t)=t^{1\over 2}/\ln{t}$ on RGG while it is
$t^{1\over 2}$ on regular lattice. More generally,
concerning the spatial propagation of opinions in social systems
or agreement dynamics in networks of artificial agents, random
geometric graphs are more realistic for a number of reasons:
({\it i}) RGG is isotropic (on average) while regular lattice is not;
({\it ii}) the average degree $\langle k \rangle$ for an RGG can be
set to an arbitrary positive number, instead of a small fixed number
for the lattice; ({\it iii}) RGGs closely capture the topology of
random networks of short-range-connected spatially-embedded
artificial agents, such as sensor networks.

An important feature of the NG which makes it different from other
models of opinion dynamics, e.g., voter model, is the spontaneous emergence of
clusters sharing the same opinion. Generally these opinion clusters
are communities closely connected within the
network. This feature of the NG can be used to detect communities of
social networks \cite{Lu1}. For the NG on random geometric graph
concerned in this paper, the clusters form a spatial structure
to which we refer to as {\it opinion domains} and which are geographic regions in which all nodes share the same
opinion. A number of relevant properties of the NG on random
geometric graph have been studied by direct individual-based
simulations and discussed in \cite{Lu2,Lu3}, such as the scaling
behavior of the consensus time and the distribution of opinion
domain size. In contrast to these previous works, here we
develop a coarse-grained approach and focus on the spatial structure
of the two-word LO-NG, such as the correlation length, shape,
and propagation of the opinion domains.

In this paper, we provide the mean-field (or coarse-grained)
equation for the NG dynamics on RGG. By analyzing the mean-field
equation, we list all possible stationary solutions and find that
the NG may get stuck in stripe-like metastable states rather than
achieve total consensus. We find significant two-time-scale
separation of the dynamics and retrieve the slow process governed
by reaction-diffusion system. Using this framework, we identify
similarities and differences between NG and other relevant models of
opinion dynamics, such as voter model, majority game and Glauber ordering.

Next, we present the governing rule of the opinion domain evolution,
that in the late stage of dynamics, the propagation speed of the
opinion domain boundary is proportional to its curvature.
Thus, an opinion domain can be considered as a mean curvature
flow making many results of the previous works applicable here
\cite{Spirin,Barros,STANLEY OSHER}. Finally we investigate the
impact of committed agents. The critical fraction of committed
agents found in the case of the NG on complete graph is also present
here. We discuss the dependence of the consensus time on the system size, the committed fraction and the average degree.


\section*{Results}

We begin with the definitions of the essential concepts of our model.

{\bf Random geometric graph} consists of $N$ agents randomly
distributed in a unit square $D=[0,1)^2$. Each agent has an
interaction range defined by $B_r(x,y)=\{(x',y')| 0<||(x'-x,
y'-y)||<r \}$, where $r$ is the local interaction radius. Two
agents are connected if they fall in each others interaction range.
The choice of network topology, denoted as $D$, impacts the boundary
conditions. Some studies, like \cite{Lu2}, choose the natural topology of the
unit square which leads to the free boundary condition. In this
paper, we assume that $D$ is a torus, imposing the
periodic boundary condition. Consequently, the opinion dynamics is free of
boundary effects until the correlation length of the opinions
grows comparable to the length scale of $D$.

{\bf Microstate} of a network is given by a spin vector
$\vec{S}=(s_1,...,s_N)$ where $s_i$ represents the opinion of the
$i^{\rm th}$ individual. In the NG, the spin value is assigned as
follows:
\begin{equation}
s_i=\left\{ \begin{array}{ll}
               1  & \mbox{for A}\\
               0 & \mbox{for AB}\\
               -1 & \mbox{for B} \end{array} \right. .
\end{equation}
The evolution of microstate is given by spin updating rules: at each
time step, two neighboring agents, a speaker $i$ and a listener $j$
are randomly selected, only the listener's state is changed (LO-NG).
The word sent by the speaker $i$ is represented by $c$, $c=1$ if the
word is A and $c=-1$ if the word is B. $c$ is a random variable depending on $s_i$. The updating rule of the NG
can be written as:
\begin{equation}
s_{j}(t+1)=\max\{\min\{s_{j}(t)+c, 1\}, -1\}.
\end{equation}

{\bf Macrostate} is given by $n_A(x,y)$, $n_B(x,y)$ and $n_{AB}(x,y)$, the concentrations
of agents at the location $(x,y)$ with opinion A, B and AB, respectively, that satisfy the normalization condition
$n_A(x,y)+n_B(x,y)+n_{AB}(x,y)=1$.
We define $s(x,y)=n_A(x,y)-n_B(x,y)$ as the local order parameter (analogous to ``magnetization"), and $\mu(x,y)$ as the local mean field
\begin{equation}
\mu(x,y)={1\over \pi r^2}\iint_{B_r(x,y)}s(x',y')dx'dy'.
\label{localMF}
\end{equation}
Finally, $f(x,y)={1\over 2}(\mu(x,y)+1)$ denotes the probability for an agent to receive a word A if it is located at $(x,y)$.

Through the geographic coarsening approach discussed in more detail in {\bf Methods}, we obtain the mean-field equation describing the
evolution of macrostate
\begin{equation}
{\partial\over \partial t}\vec{n}(x,y)=G\left(\vec{n},f\right)=
                            \left(\begin{array}{c}
                           f\\
                           1-f\end{array}\right)-\left(\begin{array}{cc}
                           1 & f\\
                           1-f & 1\end{array}\right)\vec{n} \label{eqn_macro}\;,\\
\end{equation}
while the macrostate itself is defined as
\begin{equation}
\vec{n}(x,y)=\left(\begin{array}{c}
                           n_A(x,y)\\
                           n_B(x,y)\end{array}\right)\;.
\end{equation}

\subsection*{Spatial coarsening}
\begin{figure}[!hbtp]
\begin{center}
\includegraphics[width=0.95\textwidth]{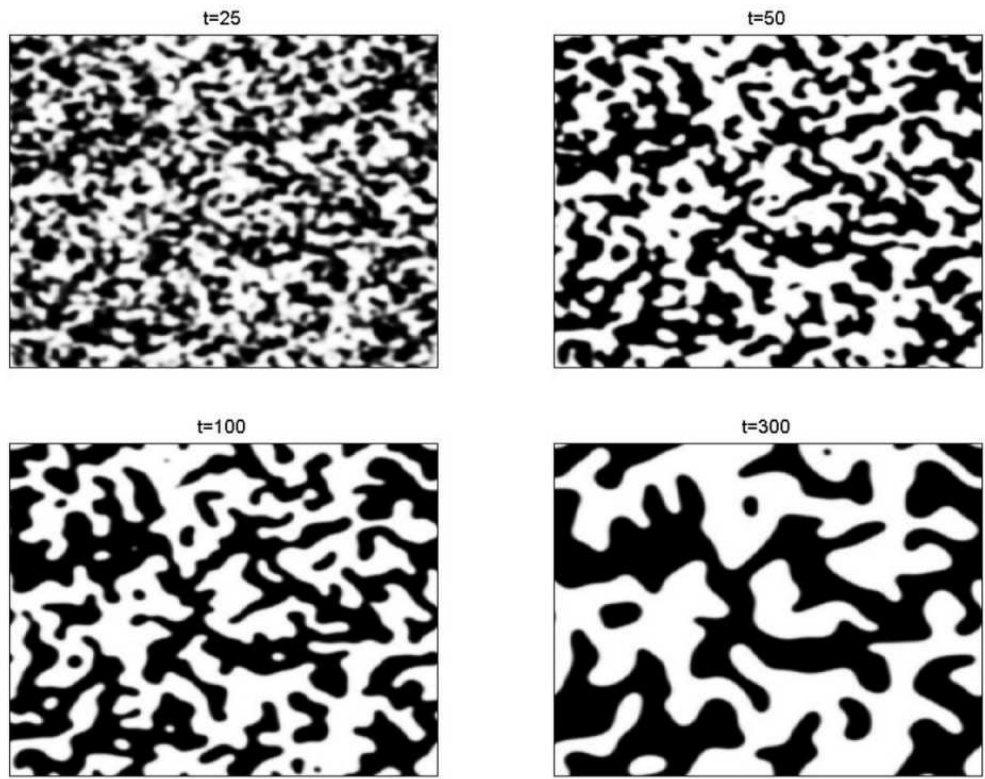}
\caption{Snapshots of numerical solution of the mean-field equation as defined by Eq.~(\ref{eqn_macro}).  Snapshots are taken at $t=25,50,100,300$, the scale of opinion domains are much bigger than $r=0.01$. Black stands for opinion A, white stands for opinion B and gray stands for the coexistence of two types of opinions. The consensus is achieved at $t\approx 10^4$.}
\label{meanfield}
\end{center}
\end{figure}

There are two characteristic length scales in this system, one is the system size (which is set to $1$), the other is the local interaction radius $r$. So regarding the correlation length or the typical scale of opinion domains $l$, the dynamics can be divided into two stages: (1) $l$ is smaller or comparable to $r$; (2) $l\gg r$.  In the second stage, the consensus is achieved when $l$ grows up to 1. Figure \ref{meanfield} present snapshots of solution of the mean-field equation. They illustrate the formation of opinion domains and the coarsening of the spatial structure.

To study the spatial coarsening, we consider the pair correlation function $C(L,t)$ defined by the conditional expectation of spin correlation.
\begin{equation}
C(L,t)=E[s(x,y,t)s(x',y',t)|\ ||(x-x',y-y')||=L].
\end{equation}

\begin{figure}[!hbt]
\begin{center}
\includegraphics[width=0.95\textwidth]{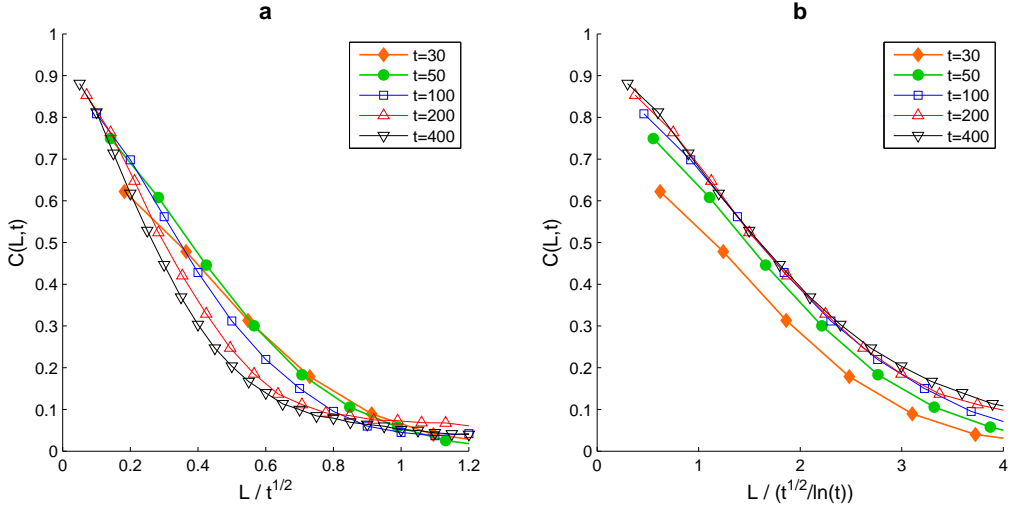}
\caption{ Scaling function $\tilde{C}(\tilde{L}=L/l(t))$ for the pair correlation function at times $t=30,50,100,200,400$. Overlapped curves indicate correct scaling of $L$. Simulations are done for the case $N=10^5$, $r=0.01$, $\langle k \rangle=31.4$. $L$ is normalized by the length scale (a) $t^{1\over 2}$ and (b)  $t^{1\over2}/\ln{t}$.}
\label{scaling}
\end{center}
\end{figure}

Figure \ref{meanfield} implies there exists a single characteristic length scale $l(t)$ so that the pair correlation function has a scaling form $C(L,t)=\tilde{C}(\tilde{L}=L/l(t))$, where the scaling function $\tilde{C}(\tilde{L})$ does not depend on time explicitly.
For coarsening in most systems with non-conserved order parameter such as the opinion dynamics on a d-dimensional lattice, the characteristic length scale is $l(t)=t^{1\over
2}$\cite{Bray}. According to the numerical results in
Fig.~\ref{scaling}, the length scale for opinion dynamics on RGG
at the early stage (t=30,50) is also $t^{1\over 2}$, but at the late
stage (t=100,200,400), the length scale $l(t)={t^{1\over 2}/
\ln{t}}$ fits more precisely simulation results than the previous one.

\subsection*{Stationary solution}
Here, we find all the possible stationary solutions of the
mean-field equation Eq.~(\ref{eqn_macro}). Taking ${\partial\over
\partial t}\vec{n}(x,y)=G(\vec{n}^*,f)=0$, we obtain
\begin{equation}
\vec{n}^*(x,y)=\left(\begin{array}{c}
                           n_A^*\\
                           n_B^*\end{array}\right)=\left(\begin{array}{c}
                            {f^2\over f^2 - f + 1}\\
                            {(1-f)^2\over f^2 - f + 1}\end{array}\right).\label{local_equil}
\end{equation}

The eigenvalues of the linear dynamical system $\lambda=- 1\pm \sqrt{f(1-f)}$ are both negative, so $\vec{n}^*$ is stable. Applying the definition of $s(x,y)$ and $\mu(x,y)$, we have
\begin{equation}
s(x,y)=n_A(x,y)-n_B(x,y)=(4\mu(x,y))/(\mu^2(x,y) + 3)
\end{equation}
\begin{equation}
\mu(x,y)={1\over \pi r^2}\iint_{ B_r(x,y)}{4\mu(x',y')\over \mu^2(x',y') + 3}dx'dy'
\end{equation}

\begin{figure}[!hbt]
\begin{center}
\includegraphics[width=\textwidth]{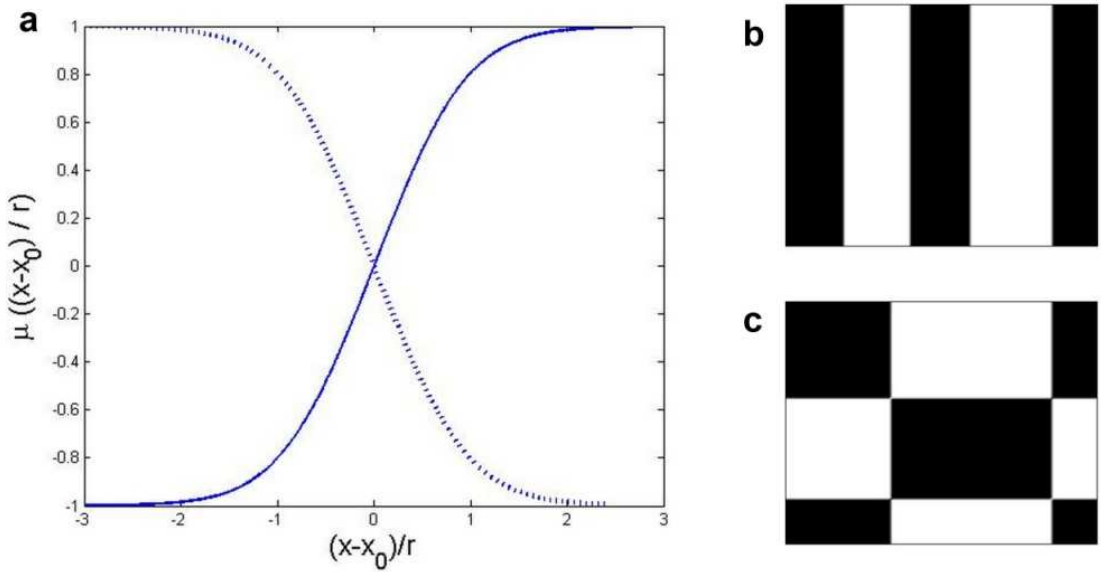}
\caption{Stationary solution. (a) Two types of intermediate layers for stationary solution $\mu({x-x_0\over r})$. $x_0$ is the location of the intermediate layer. The slope of the intermediate layer at $x=x_0$ is about $\gamma^*/r=1.034/r$. (b) Stripe-like stationary solution, neutrally stable. (c) Checker-board-like stationary solution, unstable.}
\label{stationary_solution}
\end{center}
\end{figure}

Once we solve the above integral equation, we can retrieve the
stationary macrostate $n^*(x,y)$ by Eq.~(\ref{local_equil}). Taking
$\mu(x,y)$ as a constant, we find three solutions $\mu(x,y)=\pm 1$
or $0$. $\mu(x,y)=\pm 1$ are both asymptotically stable, while
$\mu(x,y)=0$ is unstable. Another class of solution is obtained by
taking $\mu(x,y)=\mu(x)$ (or similarly $\mu(x,y)=\mu(y)$). The
solution consists of an even number of stripe-like opinion domains
demarcated by two types of straight intermediate layers parallel to
one side of the unit square $D$ as shown in
Fig.~\ref{stationary_solution}(b). With the boundary condition
$\mu(-\infty)=-1, \mu(+\infty)=1$ or vice versa, we solve the two
types of intermediate layers $\mu({x-x_0\over r})$ as shown in Fig.
\ref{stationary_solution} (a). The intermediate layers are of the
scale $r$ and can be placed at arbitrary $x_0\in[0,1)$. This class
of solution is neutrally stable.
Finally, there is another class of
solution shown in Fig.~\ref{stationary_solution}(c) with
intermediate layers both in $x$ and $y$ directions and opinion
domains assigned as a checker board. This type of solution is
unstable at the intersections of two types of intermediate layers.
 The latter two classes of solutions can be easily generalized to the cases when the intermediate layers are not parallel to x or y axis.
Later we will show that in stationary solutions all the curvature of
the opinion domain boundary has to be 0, so the solutions mentioned
above are the only possible stationary solutions.

In conclusion, considering the stability, the final state of the
macrostate dynamics can be: (1) all A or all B consensus states
which are both asymptotically stable, (2) stripe-like solution. The
probability for the dynamics stuck in the stripe-like state before
achieving full consensus is roughly ${1\over 3}$ in analogy to
similar cases in continuum percolation and spin
dynamics \cite{Spirin,Barros}.

\subsection*{Two-time-scale separation of mean-field equation and comparison with other models}

One important observation regarding the macrostate dynamics is that the
change of local mean field is usually much slower than the
convergence of local macrostate $\vec{n}$ to its local equilibrium
$\vec{n}^*$. Let $n_A=n_A^*+\delta n_A$, $n_B=n_{B}^*+\delta n_B$,
$\delta \vec{n}=(\delta n_A, \delta n_B)^T$.  The following equation
shows the exponential rate of convergence of $\vec{n}$:
\begin{equation}
{\partial \over \partial t}\delta \vec{n}(x,y)= -\left(\begin{array}{cc}
                           1 & f(x,y)\\
                           1-f(x,y) & 1\end{array}\right)\delta \vec{n}(x,y)\label{to_equil} \;.
\end{equation}

 The largest eigenvalue is $\lambda={-1+\sqrt{f(1-f)}}\leq -{1\over 2}$. So the typical time scale $\tau_n$ of the convergence of the local macrostate is $-{1\over \lambda}\leq 2$ which is independent of time and system size. The typical time scale $\tau_f$ of the change of the local mean field is inversely proportional to the propagation speed $v$ of opinion domain boundaries, and as we will show later, is of the order $O({R\over r^2})$ where $R$ is the curvature of the opinion domain boundary. Therefore, $\tau_f\gg \tau_n$ for both long time ($R$ grows to infinity along with the time $t$) and big systems (in the sense that $r\ll 1$).

\begin{figure}[!hbt]
\begin{center}
\vspace{0cm}
\includegraphics[width=0.75\textwidth]{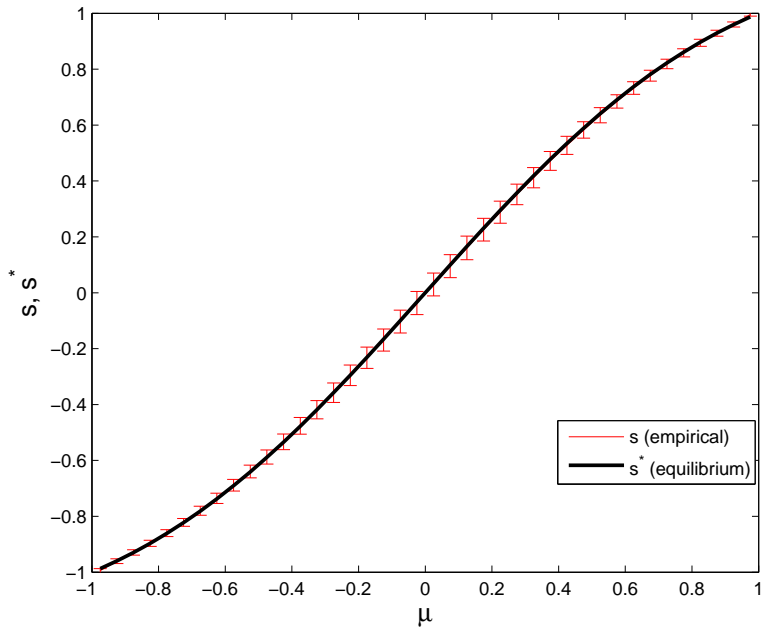}\\
\caption{Comparing the empirical local order parameter $s$ from numerical simulation and the equilibrium $s^*$ for different local mean
fields $\mu$'s. The solid line stands for $s^*$. The error bars present the means and standard deviations of $s$ that is the empirical
local order parameter for the given $\mu$ in numerical simulations.}
\label{adiabatic}
\end{center}
\end{figure}

Fig.~\ref{adiabatic} shows the significant two-time-scale separation
observed in numerical results. The equilibrium value of the local
order parameter, $s^*$, can be predicted by the local mean field,
$s^*=n_A^*-n_B^*={4 \mu/(\mu^2+3)}$. In Fig.~\ref{adiabatic}, we
present the empirical local order parameter $s$ for different local
mean field values $\mu$ and show that it is very close to its local
equilibrium $s^*$.

Since ${\partial \over \partial t}s(x,y)={\partial \over \partial t}n_A(x,y)-{\partial \over \partial t}n_B(x,y)$, we have
\begin{equation}
{\partial \over \partial t}s(x,y)={1\over 2}\left[(\mu-s) + n_{AB}\mu\right] \label{Eq_s}.
\end{equation}

This ODE is quite relevant to reaction-diffusion systems. On the right hand side, the coefficient ${1\over 2}$ is easy to get rid of
by scaling the time $\tau=t/2$. After the scaling, the first term is diffusive since $(\mu-s)$ is the continuum approximation of the
Laplace operator on RGG network acting on $s$. The second term $n_{AB}\mu$ is the local reaction term. Though classified rigorously,
it is non-local as defined in reaction-diffusion system, it represents a reaction in local neighborhood $B_r(x,y)$. The adiabaticity of
the dynamics implies that the diffusion is much slower than the local reaction. We can obtain an approximated ODE for slow time scale dynamics in a closed form by estimating $n_{AB}$ by its local equilibrium $1-n*_A-n*_B={f(1-f)\over 1-f(1-f)}={1-\mu^2\over 3+\mu^2}$.
\begin{equation}
{\partial \over \partial \tau}s(x,y)=(\mu-s) + {\mu(1-\mu^2)\over 3+\mu^2} \label{slow}.
\end{equation}

\begin{figure}[!hbt]
\begin{center}
\vspace{0cm}
\includegraphics[width=0.75\textwidth]{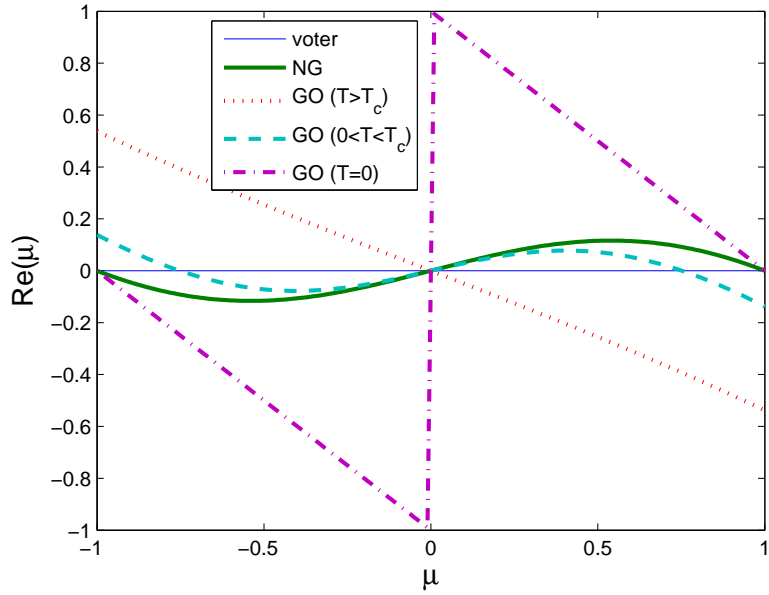}\\
\caption{Reaction term $Re(\mu)$ for voter, NG and Glauber ordering (GO) at different temperatures.}
\label{reaction_term}
\end{center}
\end{figure}

The qualitative behavior of the reaction-diffusion system is
determined by the linear stability of the reaction term
\cite{Vazquez}. In this sense, Eq.~(\ref{slow}) provides clear differentiation
between dynamics in our model and in the voter
model, the majority game and the Glauber ordering. Taking a similar
approach, we find that the voter model on RGG is purely diffusive, i.e. the
reaction term is 0. For the Glauber ordering, the reaction term is
$\tanh(\beta J\mu)-\mu$ in which $\beta$ is the inverse of temperature and
$J$ is the interaction intensity. Fig.~\ref{reaction_term} shows the
reaction term $Re(\mu)$ for the voter, NG, and Glauber ordering (GO)
at different temperatures. The majority game, NG, and Glauber ordering
at zero temperature have reaction terms with the same equilibria and
stability ($\pm 1$ stable, $0$ unstable). Thus, the mean-field solutions of
these models behave similarly. However, at the level of the
discrete model, the NG on RGG will always go to a microstate
corresponding to some stationary mean-field solution, while Glauber
ordering at zero temperature on RGG may get stuck in one of many
local minima of its Hamiltonian.

\subsection*{Boundary Evolution}
The evolution of the opinion domains is governed by a very simple
rule. The boundary of opinion domains propagate at the speed $v$
that is proportional to its curvature $1/R$, i.e. $v={\alpha\over
R}$. Here, $\alpha$ is a constant defined by the average
degree $\langle k \rangle$. In {\bf Methods}, we provide a heuristic
argument and using the perturbation method prove this
relation for the mean-field equation, i.e. the case when $\langle k
\rangle \rightarrow \infty$. This behavior is common for many
reaction-diffusion systems and it is qualitatively the same as the behavior
of Glauber ordering at zero temperature \cite{Bray,STANLEY OSHER}.

Following the rule of boundary evolution minimizes the length of the
domain boundary. A direct
consequence of this fact is that if any stationary solution exists, its
boundaries must be all straight (geodesic), confirming our conclusion about
the stationary solutions found in the previous paragraph.
Since global topology is irrelevant to our derivation,
this relation applies also to other two-dimensional manifolds.
The manifold considered here is the torus embedded in 2D Euclidean
space. However, for the standard torus embedded in
3D Euclidean space, the topology is the same but metrics are not,
hence the geodesics are different. Therefore, there are quite different and
more complicated stationary solutions there. Another example is the
sphere in 3D space. On the sphere, the only inhomogeneous stationary
solution consists of two hemispherical opinion domains, because the
great circle is the only closed geodesic on a sphere.

\begin{figure}[!hbt]
\begin{center}
\includegraphics[width=0.8\textwidth]{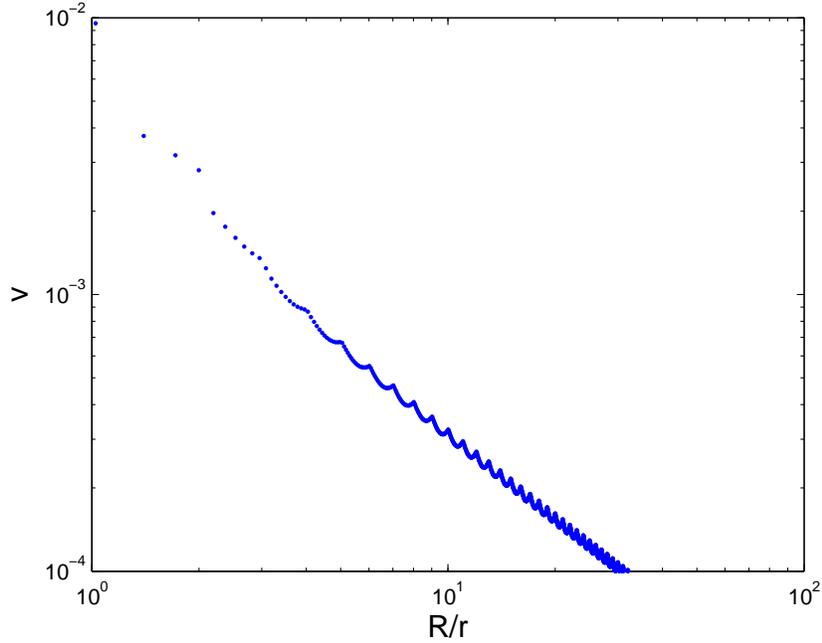}\\
\caption{Propagation speed $v$ of the boundary of opinion domains vs. curvature radius $R$. Data points are gathered from 100 runs of macrostate equation with different initial conditions. In each run, propagation speed and curvature radius are calculated for 10000 points on the boundary.}
\label{propagation_speed}
\end{center}
\end{figure}

The numerical result presented in Fig.~\ref{propagation_speed}
confirms this relation. In Fig.~\ref{propagation_speed}, we gather
$10^6$ data points from numerical solutions of the macrostate
equation, using different initial conditions, taking snapshots at
different times, tracking different points on the boundary and
calculating the local curvature radius $R$ and boundary propagation
speed $v$. These data points in the double-log plot are aligned well
with the straight line with slope $-1$. The curve formed by data points is
jiggling with some period because we implemented the numerical method
on a square lattice, so the numerical propagation speed is slightly
anisotropic.

\begin{figure}[!hbt]
\begin{center}
\includegraphics[width=0.8\textwidth]{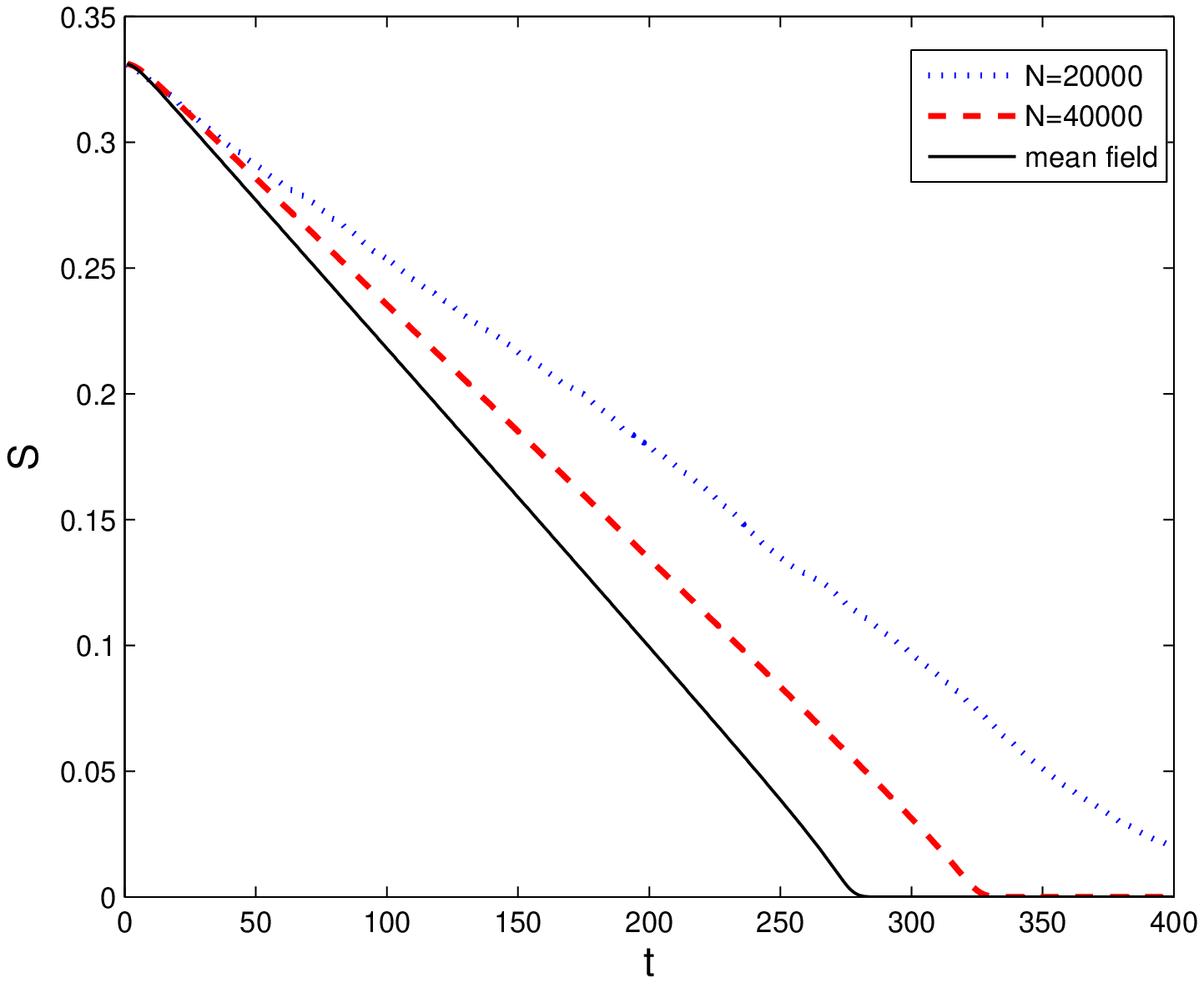}\\
\caption{The size of an opinion domain $S$ as a function of time with $r=0.05$. The initial radius of the opinion domain is $R_0=0.325$. The straight line is for the numerical solution of
mean-field equation. The dotted line and dash line are for simulation of discrete model with
$N=10000$ and $\langle k\rangle=78.5$, as well as with $N=20000$ and $\langle k \rangle=157)$, respectively.}
\label{shrink}
\end{center}
\end{figure}

Another way to confirm this rule is to consider a round opinion
domain with initial radius $R_0$. Given $v={\alpha \over R}$, the
size of this opinion domain decreases as $S(t)=\pi R^2(t)=\pi
R_0^2-2\pi \alpha t$  when $t<{R^2(0)\over 2\alpha}$. This relation
is shown in Fig.~\ref{shrink}. In simulations, the size of opinion
domain $S(t)$ is evaluated by ${1\over N}(N_A+{1\over 2}N_{AB})$
where $N_A$,$N_{AB}$ are the total numbers of A, AB nodes,
respectively. We also observe from this plot that $\alpha$ increases
with average degree $\langle k\rangle$ and converges to its upper
bound when $\langle k\rangle\rightarrow \infty$.

\subsection*{Impact of committed agents}
We now consider influencing the consensus by committed agents. In
sociological interpretation, a committed agent is one who keeps its
opinion unchanged forever regardless of its interactions with other agents.
The effect of committed agents in
the NG has been studied in \cite{Zhang1,Xie,Zhang2}. A critical
fraction of committed agent $q_c=7-4\sqrt{3}\approx 0.0718$ is found
for NG-LO on complete graph which is also relevant here.  In our
setting, a fraction $q$ of agents are committed to opinion A, and
all other agents are uncommitted and hold opinion B. The macrostate with committed
agents is still defined by Eq.~(\ref{eqn_macro}), but the definition of local
mean field $\mu$ is replaced by

\begin{equation}
\mu(x,y)={1\over \pi r^2}\iint_{ B_r(x,y)}(1-q){4\mu(x',y')\over \mu^2(x',y') + 3}+ q dx'dy'  .
\end{equation}

Generally, $q$ can vary on the x-y plane, but we only consider the
case that $q$ is a constant, i.e. the committed agents are uniformly
distributed. Then we reanalyze the stationary solution. Firstly
there is a critical committed fraction $q_c$ which is exactly the
same as that on complete graph. When $q>q_c$, the only stationary
solution is $\mu(x,y)=1$ and it is stable. When $q<q_c$, there are three solutions, of which two
$\mu(x,y)=1$, ${q-1-\sqrt{q^2-14q+1} \over 2}$ are stable, and the third
$\mu(x,y)={q-1+\sqrt{q^2-14q+1} \over 2}$ is unstable. Besides there
is a class of stationary solutions when the committed fraction is below the critical.
They are analogues of the stripe-like solutions in the non-committed agent
case. The evolution of the boundary can be interpreted as
a mean curvature flow. In such view, the fraction of agents committed to A
opinion exert a constant pressure on
the boundary surface from the side of A opinion domain. Similarly,
agents holding opinion B exert a pressure from the side of
B opinion domain. So the stationary solution will contain opinion
domains in the form of disks with critical radius $R_c$ for which the
pressure arising from the curvature offsets the pressure from the committed
fraction; thus we have ${1\over R_c}\sim q$. This type of
stationary solutions are unstable, the round disk of the opinion
domain will grow when $R>R_c$ and will shrink when $R<R_c$.
In the first case, when the typical length scale of opinion domains grows beyond
$2R_c$, the system will achieve consensus very quickly.

\begin{figure}[!hbtp]
\begin{center}
\vspace{0cm}
\includegraphics[width=0.8\textwidth]{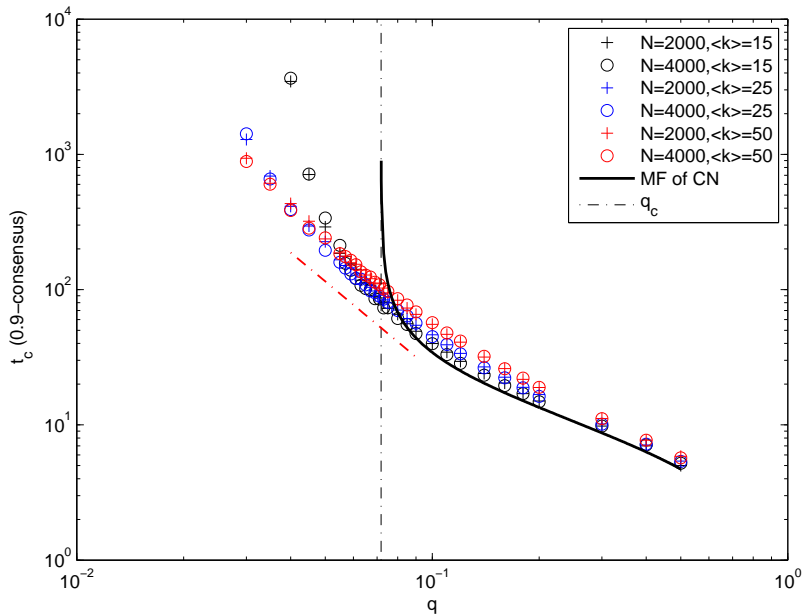}\\
\caption{Time $t_c$ for $0.9$-consensus for NG on RGG with different fractions ($q$) of committed agents with
direct simulation on networks with average degrees $\langle k\rangle=15,25,50$, and network sizes $N=2000,4000$. The solid black curve shows the mean field prediction of NG on the complete network(CN) as the limit case when $\langle k\rangle\rightarrow \infty$.
The slope of the red data points near $q_c$ is -2.19.}
\label{committed2}
\end{center}
\end{figure}

On the basis of the above stability analysis, we then analyze the dependence of the consensus time on the system size $N$, the committed fraction $q$, and the average degree $\langle k\rangle$, and show our conclusions are consistent with the numerical results in Fig.~\ref{committed2} which for a given fraction $\alpha$ ($\alpha=0.9$ in the figure) depicts the time for $\alpha$-consensus  in which at least fraction $\alpha$ of agents hold the same opinion. The time to achieve $\alpha-$consensus is independent of $N$ both according to the mean field prediction and numerical plots. When $q>q_c$ the dynamics will converge to its unique local equilibrium $\mu=1$ at all locations simultaneously. The consensus in this case is close to that on the complete network, especially when $\langle k\rangle$ tends to infinity. In the opposite case, when $q<q_c$, the process to consensus is twofold - before and after the A opinion domain achieves the critical size $2R_c$. After this criticality, the process is just the extension of the opinion domain driven by the mean field Eq.~(\ref{eqn_macro}). This stage is relatively fast and the consensus time is dominated by the duration of the other stage, the one before the criticality, in which the dynamic behavior is a joint effect of the mean field and the random fluctuation we neglect in mean field analysis. Assuming the dynamics was purely driven by the random fluctuation, the typical length scale of opinion domains would have the scaling $O(t^{1\over 2})$ at the early stage, hence the time scale of this stage would be $O(R_c^2)$, i.e. $O(q^{-2})$. If the dynamics was purely driven by the mean field, the A opinion domain would never achieve the critical size. The actual dynamics behavior is in between the two extreme cases. When $\langle k\rangle$ decreases, the fluctuation level relative to the mean field becomes higher, leading to faster consensus. In Fig.~\ref{committed2}, linear regression for the data points $0.6<q<q_c$ gives $t_c\sim q^{-\gamma}$ in which $\gamma=2.59,2.34,2.19$ for $<k>=50,25,15$ respectively, where $\gamma=2$ is the value corresponding to the purely random extreme case.

Two observations here may have meaningful sociological interpretation:

(1) When $\langle k\rangle\rightarrow \infty$ , for both $q<q_c$ and $q>q_c$, the dynamics behavior converges to that on the complete networks, though the RGG itself may be far from the complete network (with $r$ kept constant, the diameter of the RGG network is much larger than 1).

(2) When $q<q_c$, committed agents are more powerful in changing the prevailing social opinion on RGGs with low average degree. It is similar to the result of the previous study\cite{Zhang2} on the social dynamics on sparse random networks but the "more powerful" is in a different sense meaning not the smaller tipping fraction (with longer consensus time), but the shorter time to consensus.

\section*{Discussion}

On RGGs, the average degree of an agent $\langle k\rangle=\pi r^2 N$
is an important structural parameter which also strongly impacts the
local dynamic behavior. There are two critical values of $\langle k\rangle$: One is for the emergence of the giant component, $k_{c1}=4.512$ \cite{Dall}; The other one, $k_{c2}$, only applicable for finite-size networks, is for the emergence of the giant component with all vertices belonging to it.
In this paper, we only considered the case
when $\langle k\rangle$ is above the critical value $k_{c2}\sim \ln{N}$
so that the network is connected \cite{Penrose}.

We mainly focused on analyzing the mean-field equation of the NG dynamics on RGG. We
predicted a number of behaviors, including the existence of metastable
states, the two-time-scale separation, and the dependence of
the boundary propagation speed on the boundary curvature. We demonstrated
in detail that the evolution of the spatial domains for the two-word
LO-NG is governed by coarsening dynamics, similar to the broader
family of generalized voter-like models with intermediate states
\cite{Baronchelli2,Lu2,Vazquez,DallAsta_JPA2008}. However, there are
still some behaviors that cannot be explained by the mean-field
equation, such as: ({\it i}) in the large $t$ limit, the scaling of
correlation length is not $t^{1\over 2}$ as on the 2-$d$ regular
lattice, but $t^{1\over 2}/\ln{t}$; ({\it ii}) the propagation speed
increases along with the average degree $\langle k\rangle$ and its
upper-bound is the mean field prediction, i.e., the $\langle
k\rangle\rightarrow \infty$ case. So the major limitation of the
mean-field equation derived from the geographic coarsening approach
is that it neglects the fluctuation among replicas (see {\bf
Methods}) and loses the information about $\langle k\rangle$. The
dependence of the dynamics on $\langle k\rangle$ is left for further study.

\section*{Methods}
\subsection*{Geographic coarsening approach}
First, we provide the equation for the evolution of microstate
$\vec{S}=(s_i)$.  Denote the probabilities for $s_i$ taking values
$1$, $0$, $-1$ as $p_{iA}$, $p_{iAB}$ and $p_{iB}$, respectively. The
master equation for spin $s_i$ is given by
\begin{equation}
\left\{\begin{array}{l}
{d\over dt}p_{iA}=f_i (1-p_{iA}-p_{iB})-(1-f_i) p_{iA} \\
{d\over dt}p_{iB}=(1-f_i) (1-p_{iA}-p_{iB})-f_i p_{iB}
\end{array} \right.\;,
\label{eqn_micro}
\end{equation}

where $f_i={\mu_i+1\over 2}$ is the probability for the $i^{\rm th}$ agent receiving
a signal $A$, while $\mu_i$ is the local mean field defined as
the average of the neighboring spins,
\begin{equation}
\mu_i={1\over k_i}{\sum_{\{j|(x_j,y_j)\in B_r(x_i, y_j)\}} s_j},
\end{equation}
where $(x_i,y_i)$ is the coordinate of the $i^{\rm th}$ agent and
$k_i$ is the degree of the $i^{\rm th}$ agent. Master
equations for all spins together describe the evolution of
microstate.

The motivation for geographic coarsening comes from the fact that RGG is embedded in a geographic space, so we may skip the level of agents and relate the opinion states directly to the geographic coordinates. Instead of taking into account the opinion of every agent, in geographic coarsening, we consider the concentration of agents with different opinions at a specific location. $\vec{n}(x,y)$, $s(x,y)$, $f(x,y)$ are continuously differentiable w.r.t. $x$, $y$ and $t$. We derive the equation of $\vec{n}(x,y)$ from Eq.~(\ref{eqn_micro}) by taking limits,
\begin{equation}
\left\{\begin{array}{l}{\partial\over \partial t}n_A(x,y)=f(x,y) (1-n_A(x,y)-n_B(x,y))-(1-f(x,y)) n_A(x,y)\\
{\partial\over \partial t}n_B(x,y)=(1-f(x,y)) (1-n_A(x,y)-n_B(x,y))-f(x,y) n_B(x,y)
\end{array} \right.\;.
\end{equation}
The limit is done as follows
\begin{equation}
n_A(x,y)=\lim_{\epsilon\rightarrow 0}{1\over \pi \epsilon^2 N K} E\left[ \sum_{\{i|(x_i,y_i)\in B_{\epsilon}(x, y)\}}1_{\{s_i=1\}} \right] \;.
\end{equation}
The coarsening based on the above limit is valid under either of the following
two assumptions. The first is when the RGG is very dense ($N \rightarrow \infty$) and $\epsilon^2 N$ keeps constant. The second is when we consider $K$ replicas of RGG with the same set of parameters, and the summation above is over all replicas. In addition, $\epsilon^2 K$ keeps constant. Our derivation is actually based on the second assumption. However under the first assumption, the fluctuation is vanishing, and the dynamics behavior of a single run converges to the mean field result.

\subsection*{Propagation speed of the domain boundary}
We show here the qualitative property of the boundary evolution by a
heuristic argument. We consider a solution with the form
$s(x,y)=g(\tilde{R})$ where $\tilde{R}=||(x-x_0,y-y_0)||$ and with
boundary conditions $g(0)=-1$ and $g(\infty)=1$. $g(\tilde{R})$ has
an intermediate layer at $R$ as shown in
Fig.~\ref{stationary_solution}(a), so near $R$, $g(\tilde{R})\approx
{\gamma^* (\tilde{R}-R)\over r}$. When $R\gg r$, using moving
coordinate $\xi=\hat{k}\cdot \vec{x}-vt$ where $\hat{k}$ is unit
wave vector, $\vec{x}$ is spatial coordinate and $v$ is the wave
speed, Eq.~(\ref{Eq_s}) becomes
\begin{equation}
(1+n_{AB})\mu-s+2v{\partial \over \partial \xi}s=0\label{traveling_wave}\ .
\end{equation}
Here, $\mu$ can be approximated by
\begin{equation}
\mu\approx {1\over \pi r^2}\int_0^{2\pi}\int_0^r {\gamma^*\over r} \left(\sqrt{R^2+r'^2+2Rr'\cos{\theta}}-R\right) r'dr'd\theta={\gamma^* r\over 6R}\ .
\end{equation}

Then we make perturbation on Eq.~(\ref{traveling_wave}) $s=s^*+\epsilon \tilde{s}$, $\mu=\mu^*+\epsilon \tilde{\mu}$, $v=v^*+\epsilon \tilde{v}$ and so on, requiring $s^*(\xi=0)=\mu^*(\xi=0)=v^*={s}(\xi=0)=0$, and obtain the equation for $O(\xi)$
\begin{equation}
(1+n_{AB}^*)\tilde{\mu}+\tilde{n_{AB}}\mu^*-\tilde{s}+2\tilde{v}{\partial \over \partial \xi}s^*\ .
\end{equation}
At $\xi=0$, $n_{AB}^*=1/3$, $\tilde{\mu}=\mu/\epsilon$, $\tilde{v}=v/\epsilon$ and ${\partial \over \partial \xi}s^*$ is $\gamma^*/r$, the above equation becomes ${2\over 3}\mu+v\gamma^*/r =0$, so
\begin{equation}
v\approx-{2\mu r\over 3 \gamma^*}\approx -{r^2 \over 9 R}\ .
\end{equation}

\section*{Acknowledgment}

This work was supported in part by the Army Research Laboratory under Cooperative Agreement Number W911NF-09-2-0053, by the Army
Research Office Grant Nos  W911NF-09-1-0254 and W911NF-12-1-0546, and by the Office of Naval Research Grant No. N00014-09-1-0607.
The views and conclusions contained in this document are those of the authors and should not be interpreted as representing the
official policies either expressed or implied of the Army Research Laboratory or the U.S. Government.

\section*{Author Contributions}
Z.W., C.L., G.K. and B.K.S. designed the research;
Z.W implemented and performed numerical experiments and simulations;
Z.W., C.L., G.K. and B.K.S. analyzed data and discussed results;
Z.W., C.L., G.K. and B.K.S. wrote and reviewed the manuscript.

\section*{Competing financial interests}
The authors declare no competing financial interests.

\end{document}